# EFFICIENT MEASUREMENT OF THE VIBRATIONAL ROGUE WAVES BY COMPRESSIVE SAMPLING BASED WAVELET ANALYSIS

CİHAN BAYINDIR

Engineering Faculty, Işık University, 34980, Şile, İstanbul, Turkey,
cihanbayindir@gmail.com

**Abstract**

Nonlinear rogue waves appear as a result of spectral stochastic supercontinuum generation in nonlinear dynamical models, such as the nonlinear Schrödinger equation (NLSE) [1]. They are observed in the fields including but are not limited to fluid mechanics and optics. Rogue wave phenomena can also be observed in solid mechanics, where the envelope of a surface vibrational wave packet of an elastic Euler-Bernoulli beam or Kirchhoff-Love plate can be modeled in the frame of the NLSE [2, 3]. Their efficient measurement and early detection is important to analyze critical displacements, stresses and resonance. One of the promising techniques for the efficient measurement and detection of those waves is to measure the triangular Fourier spectra which begin to appear at the early stages of their development [4]. For the early detection purposes it is possible to treat such a spectrum as a sparse signal due to energy located at the central wavenumbers only [5]. Therefore compressive sampling can be a very efficient tool for measuring and predicting rogue waves [5]. However, Fourier analysis can only detect whether a vibrational rogue wave will develop or not in the vibrating medium. In order to locate it is occurrence location, wavelet analysis deems necessary if a spectral analysis approach is utilized [6]. In this paper we discuss the possible usage of the compressive sampling based wavelet analysis for the efficient measurement and for the early detection of one dimensional (1D) vibrational rogue waves. We study the construction of the triangular (V-shaped) wavelet spectra using compressive samples of rogue waves that can be modeled as Peregrine and Akhmediev-Peregrine solitons. We show that triangular wavelet spectra can be sensed by compressive measurements at the early stages of the development of vibrational rogue waves. Our results may lead to development of the efficient vibrational rogue wave measurement and early sensing systems with reduced memory requirements which use the compressive sampling algorithms. In typical solid mechanics applications, compressed measurements can be acquired by randomly positioning single sensor and multisensors.

## 1. Introduction

Rogue wave studies became extensive in recent years [1, 4-11]. While this phenomena are widely studied in optics and hydrodynamics, to our best knowledge there is no rogue wave study for vibrating mechanical systems. The envelope dynamics of the elastic bodies are modeled in the frame of the celebrated NLSE [2, 3]. Therefore rational soliton solutions which are accepted as rogue wave models presented in [1] can also become evident in vibrating mechanical systems. First two order rational soliton solutions of the NLSE are known as the Peregrine soliton and the Akhmediev-Peregrine soliton [1]. Also higher order solitons do exist in the Darboux transform hierarchy [1].

Early detection and efficient measurement of such rogue waves is a very important problem in fiber optics and in hydrodynamics. It is also vital to efficiently measure and detect the emergence of such waves in vibrating mechanical systems since they may play a crucial role in large deformation, resonance and failure analysis. In this paper we focus on the efficient measurement and detection of such vibrational rogue wave phenomena.

It is known that Fourier spectral analysis can be used to detect the emergence of the rogue waves [4]. The supercontinnum generation, that is the triangularization of the spectra becomes evident before the rogue wave appears in the physical domain [4]. However the Fourier analysis can only answer if a rogue wave will develop in the measurement domain, but can not locate its emergence location due to the lacking phase information. For this purpose a wavelet analysis deems necessary if a spectral analysis approach is utilized [6]. Similar to the Fourier spectra, the wavelet spectra exhibits becomes triangular (V-shaped) before the rogue wave becomes evident in the physical domain. The advantage of the wavelet analysis over the Fourier analysis is that it can locate the emergence point of the rogue in the physical domain before it becomes evident in time [6]. On the other hand, efficient measurement such phenomena is a vital problem to reduce the memory requirements of the measurement tools and thus increase their efficieny and reduce the cost. In a recent study we showed that the compressive sampling (CS) can be a very efficient tool for the measurement of the Fourier spectra of the rogue waves since it satisfies the sparsity requirement [5]. In this study we show that the wavelet spectra of the vibrational rogue waves can be reconstructed using CS since the wavelet spectra is sparse. Thus compared to the classical measurement systems, the CS based wavelet analysis can locate the emergence point of the vibrational rogue waves very efficiently before they become evident in time. Our results can be used in creating efficient and low cost measurement and early detection of rogue wave systems, with possible applications in vibrating mechanical systems as well as in optics and hydrodynamics.

## 2. Envelope Dynamics of an Elastic Body

Envelope dynamics of an elastic body can be modeled by the nonlinear Schrödinger equation (NLSE) [2, 3]. One of the most common nondimensional forms of the NLSE is given as

$$i\psi_t + \frac{1}{2}\psi_{xx} + |\psi|^2 \psi = 0 \qquad (1)$$

where x, t are the spatial and temporal variables, i denotes the imaginary number and $\psi$ is complex amplitude which denotes the envelope of the vibrating elastic body. Integrability of the NLSE is extensively studied within last forty years and many exact solutions are derived [1]. Additionally many numerical solution techniques, such as the spectral methods, for arbitrary conditions of the NLSE are studied [12-21]. Among the analytical solutions one interesting class of solutions is the rational rogue wave solutions [1]. These rational soliton solutions are considered as accurate rogue wave models [1]. One of the most early forms of the rational soliton solution of the NLSE is the Peregrine soliton [21]. It is given by

$$\psi_1 = \left[1 - 4\frac{1+2it}{1+4x^2+4t^2}\right]\exp[it] \qquad (2)$$

where t is the time and x is the space parameter. It is shown that Peregrine soliton is only a first order rational soliton solution of the NLSE. Higher order rational solutions of the NLSE and a hierarchy of obtaining those rational solutions based on Darboux transformations are given in [1].

It is known that rogue waves in a wave field may be in the form of first, second or higher order rational solutions of the NLSE [1, 5]. The form of the second order rogue wave, which is also known as the Akhmediev-Peregrine soliton, which satisfies the NLSE exactly is given as [1]

$$\psi_2 = \left[1 + \frac{G_2 + itH_2}{D_2}\right]\exp[it] \quad (3)$$

where $G_2 = 3/8 - 3x^2 - 2x^4 - 9t^2 - 10t^4 - 12x^2t^2$, $H_2 = 15/4 + 6x^2 - 4x^4 - 2t^2 - 4t^4 - 8x^2t^2$ and $D_2 = 1/8\left[3/4 + 9x^2 + 4x^4 + 16/3x^6 + 33t^2 + 36t^4 + 16/3t^6 - 24t^2x^2 + 16t^2x^4 + 16t^4x^2\right]$.

Currently to our best knowledge there is no experimental study on the existence of these rogue waves in vibrating mechanical systems, such as the Euler-Bernoulli beam or Kirchhoff-Love plate. However as the governing NLSE predicts, the envelope dynamics of such systems can exhibit these rogue vibrational wave phenomena.

For the measurement and early detection of these waves Fourier spectral analysis can be used waves [4]. The supercontinnum generation in the form of the triangularization of the spectra becomes evident before the rogue wave appears in the physical domain [4], however a wavelet analysis deems necessary to locate the emergence location of the vibrational rogue waves if a spectral analysis approach is utilized [6]. There are many different wavelets such as Haar, symlet, Daubechies, coiflet, biorthogonal, Meyer etc., just to name a few. Depending on the mother wavelet function it may or may not be possible to calculate the wavelet transform of the rational soliton solutions of the NLSE analytically. However for illustrative purposes we only present numerical results using Haar wavelet with scales of 1-32. Similar to the Fourier spectra, the wavelet spectra becomes triangular (V-shaped) before the rogue wave becomes evident in the physical domain. This property may be used to detect the emergence and emergence point of a vibrational rogue wave. On the other hand, efficient measurement such vibrational rogue wave phenomena is a vital problem to reduce the memory requirements of the measurement tools and thus increase their efficieny and reduce the cost. Recently, we have showed that the compressive sampling (CS) can be a very efficient tool for the measurement of the Fourier spectra of the rogue waves since their Fourier spectra satisfy the sparsity requirement [5]. In this study we also show that the wavelet spectra of the vibrational rogue waves can be reconstructed using CS since their wavelet spectra is also sparse. Thus compared to the classical measurement systems, the CS based wavelet analysis can locate the emergence point of the vibrational rogue waves very efficiently before the rogue waves become evident in time. Our results can lead to creation and implementation of efficient and low cost measurement and early detection of rogue waves, not necessarily limited to vibrating mechanical systems but also for rogue wave phenomena such as those observed in optics and hydrodynamics.

## 3. Efficient Measurement of The Vibrational Rogue Waves by Compressive Sampling Based Wavelet Analysis

3.1. Review of Compressive Sampling

Compressive sampling (CS) is an efficient sampling technique which exploits the sparsity of the signal. Using CS it is possible to reconstruct the signal by using far fewer samples than the number of samples that the classical Shannon-Nyquist sampling theorem states [22, 23]. CS has been extensively studied as a mathematical tool in applied sciences and engineering and we only try to give a very brief summary of the CS in this section and refer the reader to [22, 23] for a comprehensive discussion and derivation of the CS.

Let η be a K-sparse signal of length N which means only K out of N elements of the signal are nonzero. η can be represented in an orthonormal basis functions with transformation matrix λ. Common transformation used in literature are Fourier, discrete cosine or wavelet transforms. Therefore it is possible to write $\eta = \lambda \hat{\eta}$ where $\hat{\eta}$ is the transformation coefficient vector. Since η is a K-sparse signal we can discard the zero coefficients and obtain $\eta_s = \lambda \hat{\eta}_s$. Here, $\eta_s$ is the signal with non-zero elements only. The underlying algorithm of the CS is that a K-sparse signal η of length N can exactly be reconstructed from M ≥ Cµ²(Φ, λ)K log(N/K) measurements with an overwhelmingly high probability. Here, C is a positive constant and µ²(Φ, λ) is coherence between the sensing basis Φ and transform basis λ [23]. Taking M random projections by using the sensing matrix Φ we obtain g = Φη. Therefore the minimization problem can be summarized as

$$\min \|\hat{\eta}\|_{l_1} \text{ under constraint } g = \phi \lambda \hat{\eta} \quad (4)$$

where $\|\hat{\eta}\|_{l_1} = \sum_i |\hat{\eta}_i|$. So that among all signals which satisfy the given constraints, the $l_1$ solution of the CS problem becomes $\eta_{CS} = \lambda \hat{\eta}$. A more detailed discussion of the CS can be seen in [23]. It is useful to note that we are using the sparsity property of the Haar wavelet transform of ψ. Therefore we can write λψ = η where η is sparse triangular spectra, λ is the wavelet transformation matrix and $\Psi_{CS} = \hat{\eta}$ is the envelope measurement.

### 3.2. Proposed Compressive Sampling Based Wavelet Analysis for The Early Detection of Vibrational Rogue Waves

We propose to use the sparsity property of the wavelet spectra of the vibrational rogue waves in the form of the Peregrine and Akhmediev-Peregrine solitons for their efficient measurement and early detection. That is, we take M random compressive samples in the time domain and apply the CS algorithm to recover the classical sampling solution with a number of samples of N. In the next section we show that this approach can measure the vibrational rogue wave accurately, especially when the rogue wave is at its peak at t=0. Repeating this procedure at various time steps, the wavelet spectra thus the vibrational rogue wave profile in the physical domain can be acquired efficiently, therefore early detection of such wave phenomena become possible.

## 4. Compressive Sensing of the Triangular Wavelet Spectra

Supercontinuum of the Fourier spectra in the form of triangular shape can be used for the early detection of the vibrational rogue waves since they become evident at the early stages of their development [4]. This is true both for the first and the second order rational soliton solutions of the NLSE and also validated for the rogue waves in a chaotic wave field in the frame of the NLSE [4]. However due to the loss of the phase information, Fourier analysis can not locate the emergence point of the vibrational rogue wave. For this purpose, wavelet analysis deems necessary when a spectral analysis approach is utilized. The Haar wavelet analysis results for the Peregrine soliton are depicted in Fig. 1, where the first subfigure shows the 1D Peregrine soliton in the physical domain at t = 0 and t = 3. The second and the third subfigures show the Haar wavelet spectra obtained at t = 0 and t = 3, respectively. Since the temporal dynamics of the Peregrine soliton is symmetric, without loss of generality it is possible to state that triangular Haar wavelet spectra can be used for the early detection of Peregrine soliton, due to the principle of time reversal.

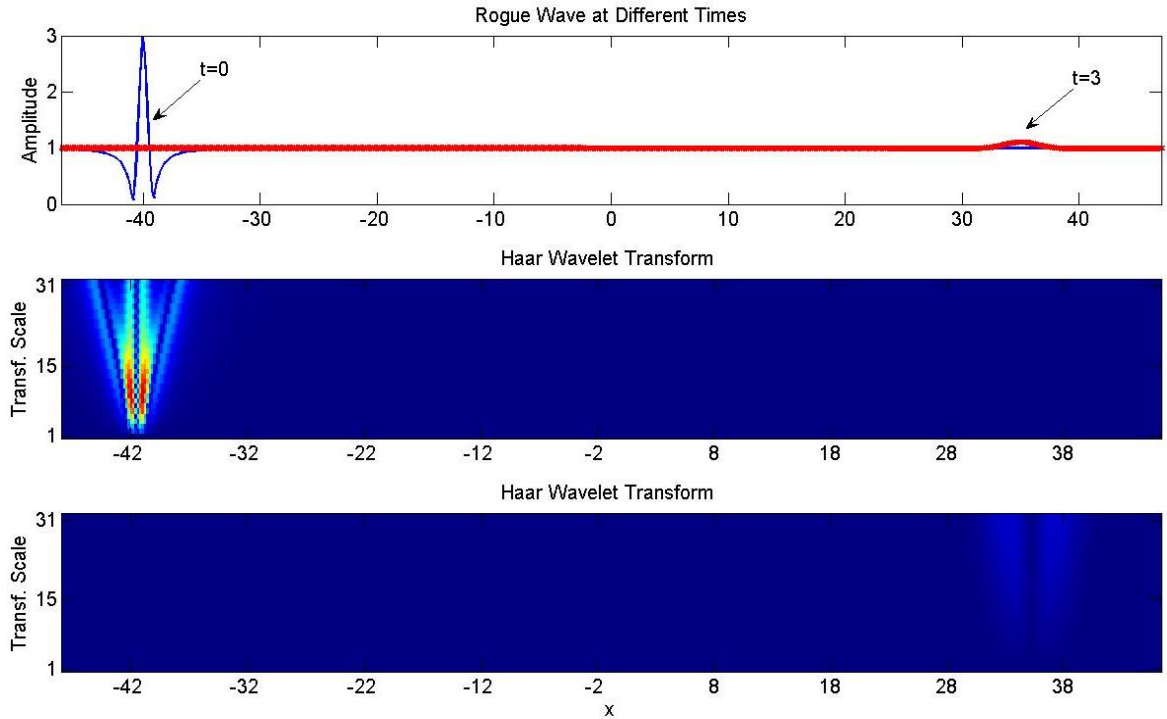

*Figure.1.* *Haar Wavelet Analysis of Peregrine Soliton at t=0 and t=3.*

For the utilization the proposed methodology, we can treat the wavelet spectra of the Peregrine soliton depicted in Fig.1 as a sparse signal. Therefore by taking random compressive measurements in physical space we can recover the sparse wavelet spectra using efficient compressive sampling algorithm. We show the results obtained this way in Fig.2., where the first subfigure shows the Peregrine soliton obtained classically using N=1024 samples and the second subfigure shows the same soliton obtained by compressive sensing using random M=64 compressive samples. The results are depicted for two different times, namely t=0 and t=3. The normalized root-mean-square (rms) difference between the two solitons in the physical domain obtained by the classical and compressive sampling techniques at time t = 0 is $9.15 \times 10^{-11}$ whereas the rms difference between two results is $7.91 \times 10^{-2}$ at time t = 3. The accuracy is reduced for the soliton at t=3 due to wider wavelet spectra which distorts the sparsity requirement. However it is still reasonably accurate and very efficient to sense the emergence of the rogue vibrational wave using the proposed technique. Additionaly the accuracy is very good at t=0, which makes the proposed technique a very efficient technique for sensing of the vibrational rogue waves. This result may be used for the development of the efficient vibrational rogue wave systems to analyze large deformations, resonance and failure of the mechanical systems. In our simulations we observe that CS can be used to obtain the sparse spectra more efficiently by employing fewer number of compressive samples (M), when the rogue wave is at its peak.

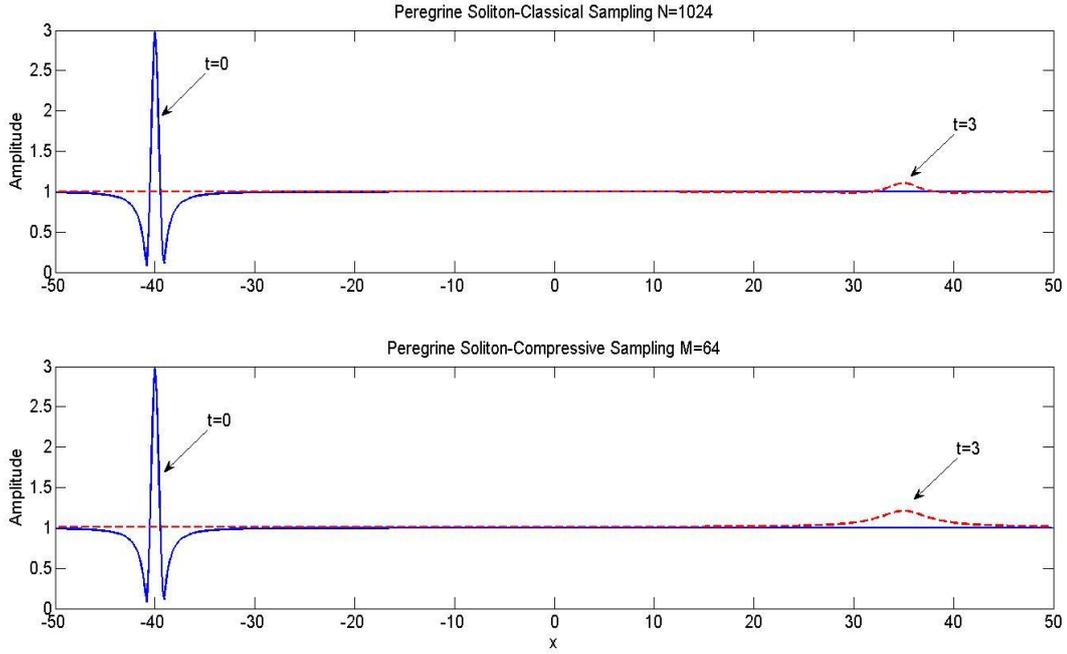

**Figure.2.** *Peregrine soliton at t=0 and t=3. a) Classical Sampling Reconstruction with N=1024 samples. b) Compressive Sampling Reconstruction with M=64 samples.*

The results for the second order rational soliton solution of the NLSE known as Akhmediev-Peregrine (AP) soliton are depicted in Fig. 3 where first subfigure shows the AP soliton in the physical domain at t = 0 and t = 3. The second and the third subfigures show its Haar wavelet spectra obtained at t = 0 and t = 3, respectively. As depicted in the figure, triangularization of the wavelet spectra of the AP soliton begins to appear at the early stages of its development. These spectra can be treated as sparse signals since the majority of their components are zero. Therefore the proposed CS based wavelet analysis can be very beneficial for measuring and early detecting the vibrational rogue waves. By taking random compressive measurements in physical space we can recover the sparse wavelet spectra using efficient compressive sampling algorithm.

We show the results for the AP soliton obtained by the proposed methodology in Fig.4., where the first subfigure shows the AP soliton obtained classically using N=1024 samples and the second subfigure shows the same soliton obtained by CS using random M=64 compressive samples. Again, the results are depicted for two different times, namely t=0 and t=3. The normalized rms difference between the two solitons in the physical domain obtained by the classical and compressive sampling techniques at time t = 0 is $8.77 \times 10^{-10}$ whereas the rms difference between two results is $9.83 \times 10^{-2}$ at time t = 3. The accuracy is again reduced at t=3 for the AP soliton due to its wider wavelet spectra which distorts the sparsity requirement of the CS algorithm. However it is still reasonably accurate and very efficient to sense the emergence of the rogue vibrational wave using the proposed technique. Additionaly the accuracy is very good at t=0, which makes the proposed technique a very effective solution for the efficient sensing of the vibrational rogue waves. We also observe that it is possible to reconstruct the sparse wavelet spectra of the AP soliton at t=0 using a smaller number of compressive samples.

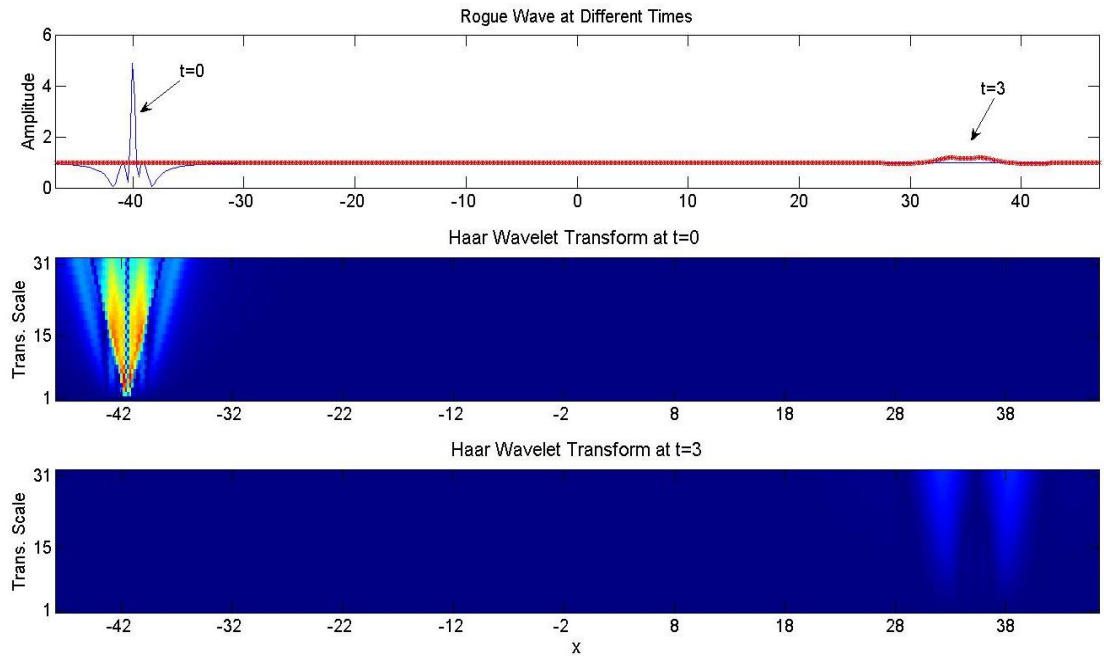

***Figure.3.*** *Haar Wavelet Analysis of Akhmediev-Peregrine Soliton at t=0 and t=3.*

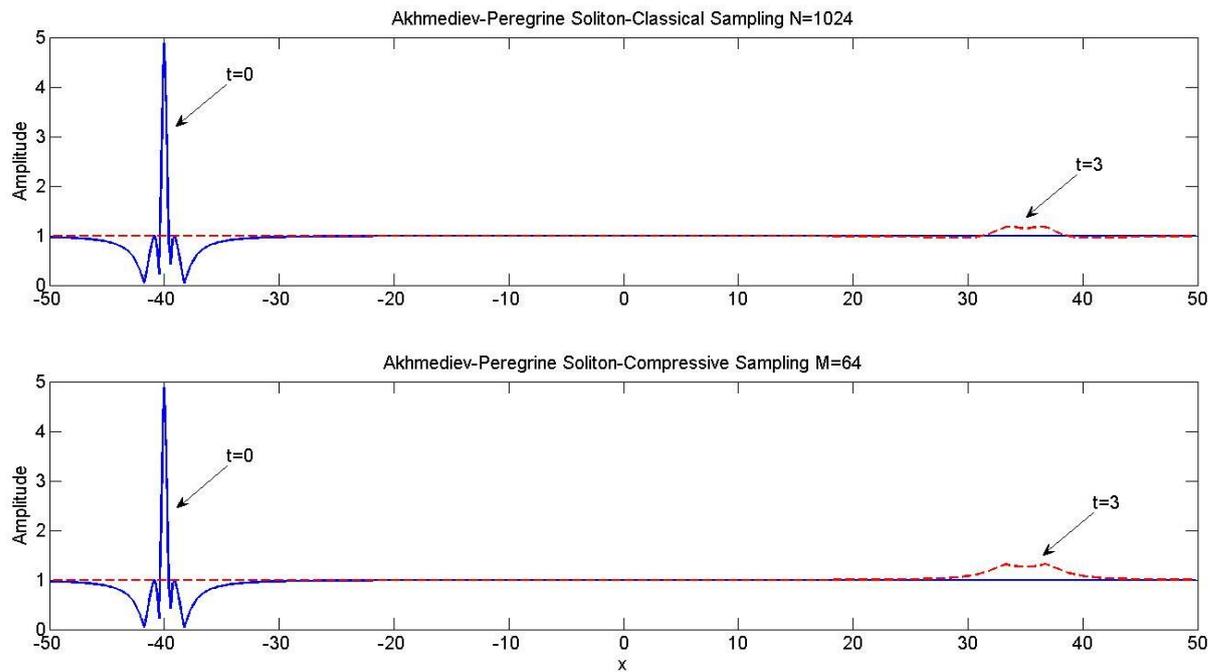

***Figure.4.*** *Akhmediev-Peregrine soliton at t=0 and t=3. a) Classical Sampling Reconstruction with N=1024 samples. b) Compressive Sampling Reconstruction with M=64 samples.*

The wavelet based CS technique proposed in this paper can be used to measure and detect the vibrational rogue waves. It may be used for the development of the efficient vibrational rogue wave systems to analyze large deformations, resonance and failure of the mechanical systems which can lead to development of the low cost sensing devices with significantly less memory requirements compared to the classical sensing ones. In practice random compressive samples can be acquired by randomly positioning single or multi sensors. Our result may also be applied in studies in the fields of fiber optics and hydrodynamics where the rogue wave phenomena are evident.

**Conclusions**

In this paper we have discussed the possible usage of the compressive sampling based wavelet analysis for the efficient measurement and the early detection of the 1D vibrational rogue waves. Particularly, we have studied the construction of the triangular wavelet spectra using compressive samples of rogue waves for which the Peregrine and Akhmediev-Peregrine solitons are accepted as accurate models. We have showed that the triangular wavelet spectra can be sensed by compressive measurements at the early stages of the development of vibrational rogue waves. Our results may lead to development of the efficient vibrational rogue wave early sensing systems with reduced memory requirements which base on the compressive sampling algorithms. In typical solid mechanics applications and experiments, the compressive measurements can be acquired by randomly positioning single or multi sensors. This technique can reduce the memory requirements of the vibration sensing systems dramatically thus will lead to development of efficient vibration sensing and early detection systems.